\documentclass[12pt, a4paper]{amsart}
\newtheorem{prop}{Proposition}[section]
\newtheorem{teo}{Theorem}[section]
\newtheorem{lema}{Lemma}[section]
\newtheorem{coro}{Corollary}[section]
\newtheorem{defi}{Definition}[section]
\newtheorem{rem}{Remark}[section]

\def\ep{\varepsilon}

\def\ve{\varepsilon}
\def\RR{{\mathbb{R}}}

\def\Oc{\RR^N\setminus\Omega}
\def\di{\displaystyle}

\begin{document}

\title[Boundary fluxes for non-local diffusion]{\bf Boundary fluxes for
non-local diffusion}

\author[C. Cortazar\and M. Elgueta \and J.D. Rossi \and N. Wolanski]{Carmen Cortazar
\and Manuel Elgueta \and Julio D. Rossi \and Noemi Wolanski}
\thanks{Supported by Universidad de Buenos Aires under grants X052 and X066,
by ANPCyT PICT No. 03-13719, Fundaci\'on Antorchas Project
13900-5, by CONICET (Argentina) and by FONDECYT (Chile).
}
\address{Carmen Cortazar and Manuel Elgueta \hfill\break\indent
Departamento  de Matem\'atica, Universidad Cat\'olica de Chile,
\hfill\break\indent Casilla 306, Correo 22, Santiago, Chile. }
\email{\tt ccortaza@mat.puc.cl, melgueta@mat.puc.cl.}

\address{Julio D. Rossi \hfill\break\indent
Consejo Superior de Investigaciones Cient\'{\i}ficas (CSIC),
\hfill\break\indent Serrano 123, Madrid, Spain,
\hfill\break\indent on leave from Departamento de Matem\'atica,
FCEyN \hfill\break\indent UBA (1428) Buenos Aires, Argentina. }
\email{{\tt jrossi@dm.uba.ar}}

\address{Noemi Wolanski \hfill\break\indent
Departamento  de Matem\'atica, FCEyN \hfill\break\indent UBA
(1428) Buenos Aires, Argentina.} \email{{\tt wolanski@dm.uba.ar} }

\date{}

\begin{abstract} We study a nonlocal diffusion operator in a bounded smooth domain
prescribing the flux through the boundary. This problem may be
seen as a generalization of the usual Neumann problem for the heat
equation. First, we prove existence, uniqueness and a comparison
principle. Next, we study the behavior of solutions for some
prescribed boundary data including blowing up ones. Finally, we
look at a nonlinear flux boundary condition.
\end{abstract}

\maketitle

\section{ Introduction}
\setcounter{equation}{0}

The purpose of this article is to address the Neumann boundary
value problem for a nonlocal diffusion equation.

Let $J: \RR^N \to \RR$ be a nonnegative, symmetric $J(z)=J(-z)$
with $\int_{\RR^N} J(z)\, dz =1$. Assume also that $J$ is strictly
positive in $B(0,d)$ and vanishes in $\RR^N \setminus B(0,d)$.
Equations of the form
\begin{equation} \label{11}
u_t (x,t) = (J*u-u) (x,t) = \int_{\RR^N} J(x-y)u(y,t) dy - u(x,t)
,
\end{equation}
and variations of it, have been recently widely used to model
diffusion processes, see \cite{BFRW}, \cite{BH}, \cite{BH2},
\cite{C},
 \cite{F},
\cite{LW},  \cite{W}. As stated in \cite{F} if $u(x,t)$ is thought
of as a density at the point $x$ at time $t$ and $J(x-y)$ is
thought of as the probability distribution of jumping from
location $y$ to location $x$, then $\int_{\RR^N} J(y-x)u(y,t)\, dy
= (J*u)(x,t)$ is the rate at which individuals are arriving at
position $x$ from all other places and $-u(x,t) = -\int_{\RR^N}
J(y-x)u(x,t)\, dy$ is the rate at which they are leaving location
$x$ to travel to all other sites. This consideration, in the
absence of external or internal sources, leads immediately to  the
fact that the density $u$ satisfies equation (\ref{11}).

Equation (\ref{11}), so called nonlocal diffusion equation, shares
many properties with the classical heat equation $u_t=\Delta u$
such as: bounded stationary solutions are constant, a maximum
principle holds for both of them and, even if $J$ is compactly
supported, perturbations propagate with infinite speed.

Given a bounded, connected and smooth domain, $\Omega$, one of the
most common boundary conditions that has been imposed to the heat
equation in the literature is the {\it Neumann boundary
condition}, $\partial u /
\partial \eta (x,t)=g(x,t)$, $x\in
\partial \Omega$.

Let us state our model  equation. We study
\begin{equation} \label{Neumann}
u_t (x,t) =  \displaystyle \int_{\Omega} J(x-y)(u(y,t) - u(x,t))
\, dy   + \displaystyle \int_{\RR^N \setminus \Omega} J(x-y)
g(y,t) \, dy ,
\end{equation}
for $x\in \Omega$. In this model we have that the first integral
takes into account the diffusion inside $\Omega$. In fact, as we
have explained the integral $\int J(x-y) (u(y,t) - u(x,t)) \, dy$
takes into account the individuals arriving or leaving position
$x$ from other places. Since we are integrating in $\Omega$, we
are imposing that diffusion takes place only in $\Omega$. The last
term takes into account the prescribed flux (given by the data
$g(x,t)$) of individuals from outside (that is individuals that
enter or leave the domain according to the sign of $g$). This is
what is called Neumann boundary conditions.

Our first result for this problem is the existence and uniqueness
of solutions and a comparison principle.

\begin{teo} \label{teo.Neumann.Lineal}
For every $u_0 \in L^1(\Omega)$ and $g\in L^\infty_{loc}
((0,\infty);L^1(\Oc))$ there exists a unique solution $u$ of
\eqref{Neumann} such that $u \in C ([0, \infty ); L^1(\Omega))$
and $u(x,0)=u_0(x)$.

Moreover the solutions satisfy the following comparison property:
$$\mbox{if }u_0(x)\le v_0(x) \mbox{ in }\Omega,  \mbox{ then }u(x,t)\le v(x,t) \mbox{
in } \Omega \times [0,\infty).$$

In addition the total mass in $\Omega$ satisfies
\begin{equation}\label{masa.Neumann}
\int_{\Omega} u(y,t) \, dy = \int_{\Omega} u_0(y) \, dy + \int_0^t
\int_\Omega \int_{\RR^N \setminus \Omega} J(x-y) g(y,s) \, dy \,
dx\, ds.
\end{equation}
\end{teo}

Once existence and uniqueness of solutions is proved an important
aspect in evolution equations is the asymptotic behavior as time
evolves. In this context, we study the asymptotic behavior of
solutions for certain fluxes on the boundary.

First, we deal with a flux independent of time, that is,
$g(x,t)=h(x)$. As happens for the heat equation, in this problem,
when $h$ verifies a compatibility condition, we prove that
solutions converge exponentially fast as $t\to \infty$ to the
unique stationary solution of the problem with the same total mass
as $u_0$. If the compatibility condition is violated then
solutions become unbounded as $t \to \infty$. We have the
following result.
\begin{teo} \label{teo.Neumann.esta}Let in addition $J\in
L^2(\RR^N)$. Let $h \in L^1(\Oc)$ such that
\begin{equation} \label{neccond.intr}
0 =  \int_\Omega  \int_{\Oc} J(x-y) h(y) \, dy \, dx.
\end{equation}
Then there exists a unique solution $\varphi$ of the problem
\begin{equation} \label{Neumann.esta}
0 =  \int_{\Omega} J(x-y)(\varphi (y) - \varphi (x)) \, dy +
\int_{\RR^N \setminus \Omega} J(x-y) h(y) \, dy
\end{equation}
that verifies $\int_\Omega u_0 = \int_\Omega \varphi$ and the
asymptotic behavior of solutions of \eqref{Neumann} is described
as follows: there exists $\beta=\beta (J, \Omega)
>0$ such that
\begin{equation}\label{conv.exp.intr}
\| u (t) - \varphi \|_{L^2 (\Omega)} \le e^{-\beta t} \| u_0 -
\varphi \|_{L^2 (\Omega)}.
\end{equation}

If \eqref{neccond.intr} does not hold then solutions of
\eqref{Neumann} are unbounded.
\end{teo}

Next, we prescribe the boundary flux in such a way that it blows
up in finite time. We consider a flux of the form
\begin{equation}\label{flujo.T-t}
g(x,t) = h(x) (T-t)^{-\alpha},
\end{equation}
with a nonnegative and nontrivial function $h$.

For this problem we analyze the possibility that the solution
becomes unbounded at time $T$ a phenomenon that is known as
blow-up in the literature. For blowing-up solutions we also
analyze the rate of blow-up (that is the speed at which solutions
go to infinity at time $T$) and the blow-up set (that is the
spatial location of the singularities).

We  find that blow-up takes place in strips of width $d $ (recall
that $J$ is positive in $B(0,d)$ and zero outside) around the
support of $h$ with blow-up rates that increase as the strips get
closer to the support of $h$.

Before stating our theorem we need some notation. We set $\Omega_0
= \Omega$, ${\mathcal B}_0 = \mbox{supp}(h)$ and define
recursively for $i \geq 1$
\begin{equation}
\begin{array}{l}
{\mathcal B}_i = \{ x\in \Omega \setminus \cup_{j<i}{\mathcal B}_j
\, : \, d(x, {\mathcal B}_{i-1} ) < d \} \nonumber
\end{array}
\end{equation}
and
$$
\Omega_i = \Omega_{i-1} \setminus {\mathcal B}_i.$$ We also define
the functions $w_i, \; \tilde w_i : \RR^N \to \mathbb{R}$
by
$$
\begin{array}{ll}
w_1 (x) =\displaystyle  \frac{1}{(\alpha -1)} \int_{\RR^N
\setminus \Omega} J(x-y) h(y) \, dy, & \\[10pt]
w_i (x) = \displaystyle \frac{1}{(\alpha -i)} \int_{\RR^N
\setminus \Omega_i} J(x-y) w_{i-1}(y) \, dy & \quad\mbox{\rm for
}1<i<\alpha, \\[10pt]
\tilde w_1 (x) =\displaystyle  \int_{\RR^N \setminus \Omega}
J(x-y) h(y) \, dy
\end{array}
$$
and
$$\tilde w_1 (x) = (\alpha -1)w_i(x) \quad\mbox{\rm for }1<i<\alpha.$$

We can now state  our result.

\begin{teo} \label{teo.T-t} Let in addition $J\in L^\infty(\RR^N)$.
Assume $h\in L^\infty(\RR^N\setminus\Omega)$, $h\ge0$,
$\int_\Omega\int_{\RR^N\setminus\Omega}
J(x-y)h(y)\,dy\,dx\neq0$. Then, the solution of \eqref{Neumann} with $g(x,t) = h(x)
(T-t)^{-\alpha}$ blows up at time $T$ if and only if $\alpha \ge
1$.

If $\alpha > 1$ is not an integer the blow-up set, $B(u)$, is
given by
$$
B(u) = \bigcup_{1\le i \leq [\alpha ]} {\mathcal B}_i,
$$
with the asymptotic behavior
$$
(T-t)^{\alpha -i } u(x,t) \to w_i (x) \; \mbox{ uniformly in }
{\mathcal B}_i \mbox{ as } \; t \to T$$ for each $i$ such that $1
\leq i \leq [\alpha ]$.

If $\alpha$ is an integer the blow-up set, $B(u)$, is given by
$$
B(u) = \bigcup_{1\le i \le \alpha} {\mathcal B}_i,
$$
with the asymptotic behavior,
$$
(T-t)^{\alpha -i } u(x,t) \to w_i (x) \; \mbox{ uniformly in }
{\mathcal B}_i \mbox{ as } \; t \to T$$ for each $i$ such that $1
\leq i < \alpha$ and
$$
\frac{u(x,t)}{-\ln (T-t)} \to \tilde w_\alpha (x) \; \mbox{
uniformly in } {\mathcal B}_\alpha \mbox{ as } \; t \to T.$$
\end{teo}

Observe that blow-up in the whole domain (global blow-up) is
possible. Indeed this happens for large values of $\alpha$
(depending on $\Omega$, $h$ and $d$).

One can compare this result with the corresponding one for the
heat equation with boundary flux $\partial u/\partial \eta (x,t) =
h(x) (T-t)^\alpha$. For the heat equation solutions blow up if and
only if $\alpha > 1/2$ and in this case $\max_{x} u(x,t) \sim
(T-t)^{-\alpha + 1/2}$. Therefore the occurrence of blow-up and
the blow-up rate for non-local diffusion are different from the
corresponding ones for the heat equation.

Finally we consider a nonlinear boundary condition of the form
\begin{equation}\label{flujo.u^p}
g(y,t) = \overline{u}^p (y,t)
\end{equation}
where $\overline{u}$ is the extension of $u$ from the boundary to
the exterior of the domain in the following form: let us assume
that a neighborhood of width $d$ of $\partial \Omega$ in $\RR^N
\setminus \Omega$ can be described by coordinates $(z,s)$ where $z
\in \partial \Omega$ and $s$ is the distance from the point to the
boundary, then we set $\overline{u} (z,s) = u(z)$. For this
nonlinear boundary condition with nonlocal diffusion we have the
following result.

\begin{teo} \label{teo.u^p} Let in addition $J\in C(\RR^N)$. Then,
positive solutions blow up in finite time if and only if $p>1$. As
for the blow up rate, there exist constants $c > 0$ and $C$ such
that
\begin{equation}\label{tasa.u^p}
c(T-t)^{-1/(p-1)} \leq \max_x u(x,t) \leq C(T-t)^{-1/(p-1)}.
\end{equation}
Moreover, the blow-up set is contained in a neighborhood of $
\partial \Omega$ of width $ K\, d $, where $K = [p/(p-1)] $.
\end{teo}

There is a large amount of literature dealing with blow-up for
parabolic equations and systems see for example the survey
\cite{GV}, the book \cite{SGKM} and references therein. When
blow-up is due to nonlinear boundary conditions see for example
\cite{Hu}, \cite{RR}, the surveys \cite{CF}, \cite{FF} and the
references therein. It is known that solutions of the heat
equation with a nonlinear boundary condition given by a power blow
up in finite time if and only if $p>1$, the blow-up rate is given
by $\| u(x,t)\|_{L^\infty (\Omega) } \sim (T-t)^{-1/(2(p-1))}$ and
the blow-up set is contained in $\partial \Omega$. Hence the
blow-up rate and set are different but the blow-up set contracts
to the boundary as the support of $J$ becomes smaller. Observe that
for $J$ fixed
the blow up set can be the whole domain $\Omega$ if $p$ is sufficiently close
to 1.

\medskip

{\bf Organization of the paper.} In Section \ref{sect-existencia}
we prove existence, uniqueness and the comparison principle, in
Section \ref{sect-h} we deal with the problem with $g(x,t)=h(x)$,
in Section \ref{SECT-b.u} we analyze the blow-up problem and
finally in Section \ref{sect-NBC} we study the problem with a
nonlinear boundary condition.


\section{Existence, uniqueness and a comparison
principle.}\label{sect-existencia} \setcounter{equation}{0}

In this section we prove Theorem \ref{teo.Neumann.Lineal} and give
as remarks several consequences of the proof that will be used
later in the paper.

As in \cite{cer}, existence and uniqueness will be a consequence
of Banach's fixed point theorem so we give first some
preliminaries.

Fix $t_0 >0$ and consider the Banach space
$$X_{t_0} =
C ([0,t_0]; L^1(\Omega))$$ with the norm
$$ ||| w|||=  \max\limits_{0\le t \le t_0} \| w(\cdot,
t) \|_{L^1(\Omega)}. $$

We will obtain the solution as a fixed point of the operator $T:
X_{t_0} \to X_{t_0}$ defined by
\begin{equation}\label{T}
\begin{array}{rl}
\displaystyle T_{w_0,g}(w) (x,t) =  w_0 (x) & + \displaystyle
\int_0^t \int_\Omega J\left( x-y \right) (w(y,s) - w(x,s))\, dy \,
ds \\[12pt]
& + \displaystyle \int_0^t \int_{\RR^N\setminus\Omega} J\left( x-y
\right) g(y,s) \, dy \, ds.
\end{array}
\end{equation}

The following lemma is the main ingredient in the proof of
existence.

\begin{lema}\label{contraction}
Let $w_0, z_0 \in L^1 (\Omega)$, $g, h \in L^\infty
((0,t_0);L^1(\Oc))$ and $w,z \in X_{t_0}$, then there exists a
constant $C$ depending only on $\Omega$ and $J$ such that
$$
|||T_{w_0,g}(w)-T_{z_0,h}(z)||| \leq \|w_0-z_0\|_{L^1(\Omega )} +
C t_0 \Big\{ |||w-z||| + \| g-h \|_{ L^\infty
((0,t_0);L^1(\Oc))}\Big\}.
$$
\end{lema}

\begin{proof} We have
$$\begin{array}{l} \displaystyle
 \int_\Omega |T_{w_0,g}(w)(x,t) - T_{z_0,h}(z)(x,t)|\, dx \leq
 \int_{\Omega} |w_0-z_0|(x) \, dx
  \\[12pt]
 \displaystyle \qquad +  \int_\Omega \left| \int_0^t
\int_\Omega J\left( x-y \right) \Big[ (w(y,s) - z(y,s))  - (w(x,s)
- z(x,s)) \Big]
 \, dy \, ds \right| \, dx \\[12pt]
  \displaystyle \qquad +  \int_\Omega \left| \int_0^t
\int_{\Oc} J\left( x-y \right) \big(g(y,s) - h(y,s)\big)
 \, dy \, ds \right| \, dx.
\end{array}$$
Hence
$$\begin{array}{l} \displaystyle
 \int_\Omega |T_{w_0,g}(w)(x,t) - T_{z_0,h}(z)(x,t)|\, dx \le
 \|w_0-z_0\|_{L^1 (\Omega)}\\[12pt]
 \displaystyle \qquad +   \int_0^t
\int_\Omega |(w(y,s) - z(y,s))|\, dy
 \displaystyle  +   \int_0^t
\int_\Omega |(w(x,s) - z(x,s))|\, dx \\[12pt]
  \displaystyle \qquad +   \int_0^t
\int_{\Oc} \left|g(y,s) - h(y,s) \right|
 \, dy \, ds  \, dx.
\end{array}$$
Therefore, we obtain,
$$
|||T_{w_0,g}(w)-T_{z_0,h}(z)||| \leq ||w_0-z_0||_{L^1(\Omega )} +
C t_0 \Big\{ |||w-z|||+ \| g-h \|_{ L^\infty
((0,t_0);L^1(\Oc))}\Big\},
$$
as we wanted to prove. \end{proof}

\begin{teo} \label{teo.L^1}
For every $u_0 \in L^1(\Omega)$ there exists a unique solution $u$
of \eqref{Neumann} such that $u \in C ([0, \infty );
L^1(\Omega))$. Moreover, the total mass in $\Omega$ verifies,
\begin{equation}\label{masa}
\int_{\Omega} u(y,t) \, dy = \int_{\Omega} u_0(y) \, dy + \int_0^t
\int_\Omega \int_{\RR^N \setminus \Omega} J(x-y) g(y,s) \, dy \,
dx\, ds.
\end{equation}
\end{teo}
\begin{proof} We check first that $T_{u_0,g}$ maps $X_{t_0}$ into
$X_{t_0}$. From \eqref{T} we see that for $0<t_1<t_2\le t_0$,
$$\begin{aligned}
\|T_{u_0,g}(w)(t_2)-T_{u_0,g}(w)(t_1)\|_{L^1(\Omega)}&\le
2\int_{t_1}^{t_2} \int_\Omega  |w(y,s)|\,dx\, dy \,
ds+\\
&+\int_{t_1}^{t_2}\int_{\RR^N\setminus\Omega} |g(y,s)|\,dx \, dy
\, ds.
\end{aligned}
$$
On the other hand, again from \eqref{T}
$$
\|T_{u_0,g}(w)(t)-w_0\|_{L^1(\Omega)}\le C t\big\{|||w||| +\|g\|_{
L^\infty ((0,t_0);L^1(\Oc))}\big\}.
$$
These two estimates give that $T_{u_0,g}(w)\in C([0,t_0];
L^1(\Omega))$. Hence $T_{u_0,g}$ maps $X_{t_0}$ into $X_{t_0}$.

Choose $t_0$ such that $ C t_0 <1$. Now taking $z_0 \equiv w_0
\equiv u_0$, $g\equiv h$ in Lemma \ref{contraction} we get that
$T_{u_0,g}$ is a strict contraction in $X_{t_0}$ and the existence
and uniqueness part of the theorem follows from Banach's fixed
point theorem in the interval $[0,t_0]$. To extend the solution to
$[0,\infty)$ we may take as initial data $u(x,t_0)\in L^1(\Omega)$
and obtain a solution up $[0,2t_0]$. Iterating this procedure we
get a solution defined in $[0,\infty)$.

We finally prove that if $u$ is the solution, then the integral in
$\Omega$ of $u$ satisfies \eqref{masa}. Since
$$
\begin{array}{rl}
u(x,t) - u_0 (x)= &  \displaystyle \int_0^t \int_\Omega
J\left(x-y\right) (u(y,s)
- u(x,s))\, dy \, ds \\[12pt]
& + \displaystyle \int_0^t \int_{\RR^N\setminus\Omega} J\left( x-y
\right) g(y,s) \, dy \, ds.
\end{array}
$$
We can integrate in $x$ and apply Fubini's theorem to obtain
$$
\displaystyle\int_\Omega u(x,t)dx - \int_\Omega u_0 (x)dx =
\displaystyle \int_0^t \int_\Omega \int_{\RR^N\setminus\Omega}
J\left( x-y \right) g(y,s) \, dy \, dx \, ds .
$$
and the theorem is proved. \end{proof}

Now we give some consequences that we state as remarks for the
sake of future references.

\begin{rem}\label{continuous}
Solutions of \eqref{Neumann} depend continuously on the initial
condition and boundary data. Let $u$ be a solution of
\eqref{Neumann} with initial datum $u_0$ and $v$ a solution of
\eqref{Neumann} with $g$ replaced by $h$ and initial datum $v_0$.
Then for every $t_0>0$ there exists a constant $C= C(t_0)$ such
that
$$\max\limits_{0\le t \le t_0} \| u(\cdot,
t)- v(\cdot ,t) \|_{L^1(\Omega)}
 \le C ||u(\cdot ,0)-v(\cdot, 0)||_{L^1(\Omega )} + C \| g-h \|_{ L^\infty
((0,t_0);L^1(\Oc))}.
$$
\end{rem}

\begin{rem}\label{iff}
The function $u$ is a solution of \eqref{Neumann} if and only if
\begin{equation}\label{repre.u}
\begin{array}{rl}
u(x,t) = & \displaystyle e^{-A(x)t }u_0(x) + \int_0^t\int_\Omega
e^{-A(x)(t-s)}J(x-y )
u(y,s) \, dy \, ds \\[12pt]
& + \displaystyle \int_0^t
\int_{\RR^N\setminus\Omega}e^{-A(x)(t-s)} J( x-y ) g(y,s) \, dy \,
ds,
\end{array}
\end{equation}
where $A(x)=\int_\Omega J(x-y)\,dy$.

Observe that $A(x)\ge\alpha>0$ ($x\in \overline{\Omega}$) for a
certain constant $\alpha$.
\end{rem}

\begin{rem}\label{derivative}
{}From the previous remark we get that if $u\in
L^\infty(\Omega\times(0,T))$, $u_0 \in C^k(\overline{\Omega} )$
with $0 \leq k \leq \infty$, $g\in L^\infty(\Oc\times(0,T))$ and
$J\in W^{k,1}(\RR^N)$, then
$u(\cdot ,t) \in C^k(\overline{\Omega}\times[0,T]) )$. 

On the other hand, if $J\in L^\infty(\RR^N)$, $u_0\in
L^\infty(\Omega)$ and $g\in L^1(\Oc\times(0,T))$, there holds that
$u\in L^\infty(\Omega\times(0,T))$. (See Corollary \ref{coro.2.3}
for an explicit bound in the case of continuous solutions).
\end{rem}

We now define what we understand by sub and supersolutions.

\begin{defi} A function $u\in C([0,T);L^1((\Omega))$ is a
supersolution of \eqref{Neumann} if $u(x,0) \ge u_0(x)$ and
\begin{equation}\label{super}
u_t(x,t)\ge \int_\Omega J(x-y)\big(u(y,t)-u(x,t)\big)\,dy+
\int_{\Oc}J(x-y)g(y,t)\,dy.
\end{equation}
\end{defi}

Subsolutions are defined analogously by reversing the
inequalities.

\begin{lema} \label{compar.0}
Let $u_0\in C(\overline\Omega)$, $u_0\ge0$, and $u\in
C(\overline{\Omega}\times [0,T])$  a supersolution to
\eqref{Neumann} with $g\ge0$. Then, $u\ge0$.
\end{lema}

\begin{proof}
Assume that $u(x,t)$ is negative somewhere. Let $v(x,t)=u(x,t)+\ep
t$ with $\ep$ so small such that $v$ is still negative somewhere.
Then, If we take $(x_0,t_0)$ a point where $v$ attains its
negative minimum, there holds that $t_0>0$ and
$$\begin{aligned}
&v_t(x_0,t_0)=u_t(x_0,t_0)+\ep>\int_{\Omega} J(x-y)(u(y,t_0) -
u(x_0,t_0)) \, dy\\
&\hskip1cm=\int_{\Omega} J(x-y)(v(y,t_0) - v(x_0,t_0)) \, dy\ge0
\end{aligned}
$$
 which is a contradiction. Thus, $u\ge0$.
\end{proof}

\begin{coro} Let $J\in L^\infty(\RR^N)$. Let $u_0$
and $v_0$ in $L^1(\Omega)$ with $u_0\ge v_0$ and $g,\,h\in
L^\infty((0,T);L^1(\RR^N\setminus\Omega))$ with $g\ge h$. Let $u$
be a solution of \eqref{Neumann} with $u(x,0)=u_0$ and Neumann
datum $g$ and $v$ be a solution of \eqref{Neumann} with
$v(x,0)=v_0$ and Neumann datum $h$. Then, $u\ge v$ a.e.
\end{coro}

\begin{proof}
Let $w=u-v$. Then, $w$ is a supersolution with initial datum
$u_0-v_0 \ge 0$ and boundary datum $g-h\ge 0$. Using the
continuity of solutions with respect to the initial  and Neumann
data and the fact that $J\in L^\infty(\RR^N)$, we may assume that
$u,\, v \in C(\overline{\Omega}\times
[0,T])$. By Lemma \ref{compar.0} we obtain that $w=u-v \ge 0$. So
the corollary is proved.
\end{proof}

\begin{coro} Let  $u\in C(\overline{\Omega}\times [0,T])$ (resp. $v$) be a
supersolution (resp. subsolution) of \eqref{Neumann}. Then, $u\ge
v$.
\end{coro}

\begin{proof} It follows the lines of the proof of the previous corollary.
\end{proof}

\begin{coro} \label{coro.2.3}
Let $u$ be a  continuous solution of \eqref{Neumann} with
$u(x,0)=u_0$ and Neumann datum $g\in L^\infty\big((\RR^N\setminus\Omega)\times
(0,T)\big)$. Then,
\begin{equation}\label{cota}
u(x,t)\le \sup_\Omega u_0+\int_0^t \sup_\Omega \int_{\Oc}
J(x-y)g(y,s)\,dy\,ds.
\end{equation}
\end{coro}

\begin{proof} Let
$$
v(t) = \sup_\Omega u_0+\int_0^t \sup_\Omega \int_{\Oc}
J(x-y)g(y,s)\,dy\,ds.
$$
Then,  $v$ is a continuous supersolution of \eqref{Neumann}. By
the previous corollary we get the estimate \eqref{cota}.
\end{proof}

\begin{coro} If $J\in L^\infty(\RR^N)$ and
$u\in C\big([0,T];L^1(\Omega)\big)$ is a solution
of \eqref{Neumann} with $u(x,0)\in L^\infty(\Omega)$,
$g\in L^\infty\big(0,T;L^1(\RR^N\setminus\Omega)\big)$ then,
\eqref{cota} holds.
\end{coro}

\begin{proof} Let $u_n$ be the solution of \eqref{Neumann}
with $u_n(x,0)=u_0^n(x)$ and Neumann datum $g_n\in L^\infty\big
((\RR^N\setminus\Omega)\times(0,T)\big)$ such that $g_n\to g$ in
$L^1\big ((\RR^N\setminus\Omega)\times(0,T)\big)$ and $u_0^n\to
u_0$ in $L^1\big (\Omega)$ with $\|u_0^n\|_{L^\infty(\Omega)}\le
\|u_0\|_{L^\infty(\Omega)}$.

The result follows from the application of Corollary
\ref{coro.2.3} to the functions $u_n\in C(\bar\Omega\times[0,T])$
and taking limits as $n \to \infty$.
\end{proof}


\section{Asymptotic behavior for $g(x,t) = h(x)$.}\label{sect-h}
\setcounter{equation}{0}

In this section we study the asymptotic behavior, as $t \to
\infty$, of the solutions of problem \eqref{Neumann} in the case
that the boundary data is time independent. So we will assume
throughout this section  that $g(x,t) = h(x)$ and that $J\in
L^2(\RR^N)$. We start by analyzing the corresponding stationary
problem so we consider the equation
\begin{equation} \label{Neumann.esta.2}
0 =  \int_{\Omega} J(x-y)(\varphi (y) - \varphi (x)) \, dy +
\int_{\Oc} J(x-y) h(y) \, dy .
\end{equation}

Integrating in $\Omega$, it is clear that a necessary condition
for the existence of a solution $\varphi$ is that
\begin{equation} \label{neccond}
0 =  \int_\Omega  \int_{\Oc} J(x-y) h(y) \, dy \, dx.
\end{equation}

We will prove, by means of Fredholm's alternative, that condition
\eqref{neccond} is sufficient for existence and that the solution
is unique up to an additive constant.

To do this we write \eqref{Neumann.esta.2} in the form
\begin{equation} \label{Neumann.esta.3}
 \varphi (x)-K (\varphi )(x) = b(x),
\end{equation}
where
$$
b(x) = a(x) \int_{\Oc} J(x-y) h(y) \, dy,
$$
$$
K (\varphi )(x) = a(x) \int_{\Omega} J(x-y)\varphi (y) \, dy
$$
and
$$
a(x) = \displaystyle \left( \displaystyle\int_{\Omega} J(x-y) \,
dy \right)^{-1}.
$$
We consider the measure
$$
d\mu = \frac{dx}{a(x)}
$$
and its corresponding space $L^2_\mu$ of square integrable
functions with respect to this measure.

We observe that, due to our assumptions on $J$, the operator $K$
maps $L^2_\mu$ into $L^2_\mu$ and as an operator $K: L^2_\mu \to
L^2_\mu$ is compact 
and self adjoint.

We look now at the kernel of $I - K$ in $L^2_\mu$. We will show
that this kernel consist only of constant functions. In fact, let
$\varphi \in \ker (I-K)$. Then $\varphi$ satisfies
$$
\varphi (x) = a(x) \int_{\Omega} J(x-y)\varphi (y) \, dy.
$$
In particular, since $J\in L^2(\RR^N)$, $\varphi$ is a continuous
function. Set $A=\max\limits_{x \in \overline{\Omega}}\varphi (x)$
and consider the set
$$\mathcal{A} = \{ x \in \overline{\Omega} \; | \; \varphi (x) = \mbox{A} \}.$$
The set $\mathcal{A}$ is clearly closed and non empty. We claim
that it is also open in $\overline{\Omega}$. Let $x_0 \in
\mathcal{A}$. We have then
$$
\varphi (x_0) = a(x_0) \int_{\Omega} J(x_0-y)\varphi (y) \, dy.
$$
Since $a(x_0) = (\int_{\Omega} J(x_0-y) \, dy)^{-1}$ and $\varphi
(y) \leq \varphi (x_0)$ this implies $\varphi (y) = \varphi (x_0)$
for all $y \in \Omega \cap B(x_0,d)$, and hence $\mathcal{A}$ is
open as claimed. Consequently, as $\Omega$ is connected,
$\mathcal{A}=\Omega$ and $\varphi$ is constant.

According to Fredholm's alternative problem \eqref{Neumann.esta.2}
has a solution if and only if
$$\int_\Omega b(x)\frac{dx}{a(x)} = 0$$
or equivalently
$$\int_\Omega \int_{\Oc} J(x-y)h(y)\, dy\, dx = 0.$$
We have proved
\begin{teo} Problem \eqref{Neumann.esta.2} has a solution
if and only if condition \eqref{neccond} holds. Moreover any two
solutions differ by an additive constant.
\end{teo}

We will address now the problem of the asymptotic behavior of the
solution of \eqref{Neumann}. The next proposition shows the
existence of a Liapunov functional for solutions of
\eqref{Neumann}. Its proof is a direct computation and will be
omitted.

\begin{prop} \label{Prop.funcional}
Let $u(x,t)$ be the solution of \eqref{Neumann}. Let us define
\begin{equation}\label{funcional}
\begin{array}{rl}
F(u) (t) = & \displaystyle \frac{1}{4} \int_\Omega \int_\Omega
J(x-y) (u(y,t) - u(x,t))^2 \, dy \, dx \\[15pt]
& \displaystyle - \int_\Omega \int_{\RR^N \setminus \Omega} J(x-y)
h(y) u(x,t) \, dy \, dx.
\end{array}
\end{equation}
Then
$$
\frac{\partial }{\partial t} F (u) (t) = - 2 \int_\Omega (u_t)^2
(x,t) \, dx.
$$
\end{prop}

We are now in a position to state and prove a result on the
asymptotic behavior of continuous solutions.

\begin{teo} \label{comp.asimp}
Let $u$ be a continuous solution of \eqref{Neumann} with $g(x,t) =
h(x)$ where $h$ satisfies the compatibility condition
\eqref{neccond}. Let $\varphi$ be the unique solution of
\eqref{Neumann.esta.2} such that $$\int_\Omega \varphi (x)\, dx =
\int_\Omega u_0(x)\, dx.$$ Then
\begin{equation} \label{comp.asimp.Neu}
u(x,t) \to \varphi (x) \mbox{ as } t \to \infty
\end{equation}
uniformly in $\overline{\Omega}$.

When \eqref{neccond} does not hold, solutions of \eqref{Neumann}
are unbounded.
\end{teo}

\begin{proof} Set $w(x,t) = u(x,t)-\varphi (x)$. Then $w$ satisfies
$$w_t (x,t)  = \displaystyle\int_{\Omega} J(x-y) w(y,t) \, dy   - w(x,t) \int_{\Omega}
J(x-y)\, dy $$ and $\int_\Omega w(x,t)dx \equiv 0$.

By the estimate given in Corollary \ref{coro.2.3} we have that
$\|w\|_{L^\infty({\Omega} \times [0,\infty)) }$ is bounded in
$\overline{\Omega} \times [0,\infty )$ by $\| u_0
-\varphi\|_{L^\infty (\Omega)}$.

Setting $A(x) = \int_{\Omega} J(x-y)\, dy$ and integrating, the
above equation can be written as
$$
 w(x,t) = e^{-A(x) t}w(x,0) + \int_0^t e^{-A(x)(t- s)} \int_{\Omega}
J(x-y) w(y,s) \, dy \, ds.
$$

We note that $A(x)$ is a smooth function and that there exists
$\alpha >0$ such that $A(x) \geq \alpha$ for all $x \in \overline{
\Omega}$. We observe that for $x_1,\; x_2 \in \overline{\Omega}$
one has
$$|e^{-A(x_1) t}-e^{-A(x_2) t}| \leq e^{-\alpha
t}t|A(x_1)-A(x_2)|.$$ With this inequality in mind it is not
difficult to obtain, via a triangle inequality argument,  the
estimate
$$|w(x_1,t)-w(x_2,t)| \leq D\, \big(|A(x_1)-A(x_2)| +
|w(x_1,0)-w(x_2,0)|\big)$$ where the constant $D$ is independent
of $t$. This implies that the functions $w(\cdot ,t)$ are
equicontinuous. Since they are also bounded, they are precompact
in the uniform convergence topology.

Let $t_n$ be a sequence such that $t_n \to \infty$ as $n \to
\infty$. Then the sequence $w(\cdot ,t_n)$ has a subsequence, that
we still denote by $w(\cdot ,t_n)$, that converges uniformly as $n
\to \infty$ to a continuous function $\psi$. A standard argument,
using the Liapunov functional of Proposition \ref{Prop.funcional},
proves that $\psi$ is a solution of the corresponding stationary
problem and hence $\psi$ is constant. As $\int_\Omega w(x,t)dx
\equiv 0$ this constant must be $0$. Since this holds for every
sequence $t_n$, with $t_n \to \infty$, we have proved that
$w(\cdot ,t) \to 0$ uniformly as $t \to \infty$ as we wanted to
show.

When \eqref{neccond} does not hold the equation satisfied by the
total mass, \eqref{masa}, implies that $u$ is unbounded.
\end{proof}

We end this section with a proof of the exponential rate of
convergence to steady states of solutions in $L^2$. This proof
does not use a Lyapunov argument. It is based on energy estimates.

First, we prove a Lemma that can be viewed as a Poincar\'e type
inequality for our operator.

\begin{lema} \label{Poincare}
There exists a constant $C>0$ such that for every $u\in
L^2(\Omega)$ it holds
$$
\int_\Omega (u(x) - \langle u\rangle )^2 \, dx \le C \int_\Omega
\int_\Omega J(x-y) (u(y)-u(x))^2 \, dy \, dx,
$$
where $\langle u\rangle$ is the mean value of $u$ in $\Omega$,
that is
$$
\langle u\rangle=\frac{1}{|\Omega|} \int_\Omega u(x) \, dx.
$$
\end{lema}

\begin{proof}
We can assume that $\langle u\rangle=0$. Now let us take a
partition of $\RR^N$ in non-overlapping cubes, $T_i$, of diameter
of length $h$. Using an approximation argument we can consider
functions $u$ that are constant on each of the cubes $T_i$,
$u|_{T_i} = a_i$. We will only consider cubes $T_i$ such that $T_i
\cap \Omega \neq \emptyset$. For this type of functions we have to
prove that there exists a constant $C$ independent of the
partition such that
$$
\sum_{i} |T_i| a_i^2  \le C \sum_i \sum_k \int_{T_i} \int_{T_k}
J(x-y)\, dy \, dx\, (a_i -a_k)^2 .
$$

Recall that there exist $\sigma >0$ and $r>0$ such that $J(x-y)
\ge \sigma$ for any $|x-y|<2r$.

If the centers of two cubes $T_i$, $T_k$, are at distance less
than $r$ and $h<r$ we have
$$
 \int_{T_i} \int_{T_k}
J(x-y)\, dy \, dx\, (a_i -a_k)^2 \ge \sigma |T_i| |T_k|
(a_i-a_k)^2.
$$

Given $T_i$, $T_k$ two cubes intersecting $\Omega$ there exists a
number $\ell$, depending only on $\Omega$ and $r$ but not on $h$,
such that there exist a collection of at most $\ell$, not
necessarily pairwise adjacent, cubes $T_{j_1},\dots,T_{j_\ell}$
intersecting $\Omega$ with $T_{j_1} = T_i$, $T_{j_\ell} = T_k$ and
such that the distance between the centers of $T_{j_m}$ and
$T_{j_{m+1}}$ is less than $r$. Since all the involved cubes have
the same measure, we have
\begin{equation}\label{cubos}
\begin{array}{rl}
|T_i| |T_k| (a_i-a_k)^2 & \displaystyle \le 2^\ell \left(
\sum_{m=1}^{\ell-1} |T_{j_m}| |T_{j_{m+1}}|
(a_{j_m} - a_{j_{m+1}})^2 \right) \\[12pt]
& \di \le \frac{2^\ell}{\sigma}  \sum_{m=1}^{\ell-1} \int_{T_{j_m}}
\int_{T_{j_{m+1}}} J(x-y)\, dy \, dx\, (a_{j_m} - a_{j_{m+1}})^2 .
\end{array}
\end{equation}

The intermediate cubes used in \eqref{cubos} corresponding to
each pair $T_i,\,T_k$  can be chosen in such a way that no pair
of cubes is used
more that a fixed number of times (depending only on the diameter
of $\Omega$ and $r$) when varying  the pairs $T_i,\,T_k$.
Therefore, there exists a constant $C$, depending only on $J$ and
$\Omega$ but not on $h$, such that
$$
\sum_i \sum_k |T_i| |T_k| (a_i -a_k)^2  \le C \sum_i \sum_k
\int_{T_i} \int_{T_k} J(x-y)\, dy \, dx\, (a_i -a_k)^2 .
$$

On the other hand, as we are assuming that
$$
\sum_{i} |T_i| a_i =0,
$$
we get
$$\sum_i \sum_k |T_i| |T_k| (a_i -a_k)^2 \ge
2 |\Omega| \sum_i|T_i| (a_i)^2
$$
and the result follows.
\end{proof}

Now let us take the best $J-$Poincar\'e constant that is given by
\begin{equation}
\label{mejor.c.Poincare} \beta = \inf_{u\in
L^2(\Omega)}\displaystyle\frac{ \displaystyle \int_\Omega
\int_\Omega J(x-y) (u(y)-u(x))^2 \, dy \,
dx}{\displaystyle\int_\Omega (u(x) - \langle u\rangle )^2 \, dx }.
\end{equation}
Note that by Lemma \ref{Poincare} $\beta$ is strictly positive and
depends only on $J$ and $\Omega$.

Now let us prove the exponential convergence of $u(x,t)$ to the
mean value of the initial datum when the boundary datum vanishes,
i.e., $h=0$.

\begin{teo} \label{teo.Neumann.Lineal.conv.Exp} For every $u_0\in L^2
(\Omega)$ the solution  $u(x,t)$ of \eqref{Neumann} with $h=0$,
satisfies
\begin{equation}\label{conv.exp}
\| u(\cdot,t) - \langle u_0\rangle \|^2_{L^2 (\Omega)} \le
e^{-\beta t} \| u_0 - \langle u_0\rangle \|^2_{L^2 (\Omega)}.
\end{equation}
Here $\beta$ is given by \eqref{mejor.c.Poincare}.
\end{teo}

\begin{proof} Let
$$
H(t) = \frac12 \int_\Omega (u(x,t) - \langle u_0\rangle)^2 \, dx.
$$
Differentiating with respect to $t$ and using
\eqref{mejor.c.Poincare}, recall that $\langle u\rangle=\langle
u_0\rangle$, we obtain
$$
H'(t) = - \frac12 \int_\Omega \int_\Omega  J(x-y) (u(y,t) -
u(x,t))^2 \, dy \, dx \le - \beta  \frac{1}{2} \int_\Omega
(u(x,t)-\langle u_0\rangle)^2 \, dx.
$$
Hence
$$
H' (t) \le - \beta H(t).
$$
Therefore, integrating we obtain,
$$
H(t) \le  e^{-\beta t} H(0).
$$
As we wanted to prove.
\end{proof}

As a corollary we obtain exponential decay to the steady state for
solutions of \eqref{Neumann} with $h\neq 0$.

\begin{coro} \label{decaimiento.h.neq.0}
For every $u_0\in L^2 (\Omega)$ the solution of \eqref{Neumann},
$u(x,t)$, verifies
\begin{equation}\label{conv.exp.coro}
\| u - \varphi \|^2_{L^2 (\Omega)} \le e^{-\beta t} \| u_0 -
\varphi \|^2_{L^2 (\Omega)} .
\end{equation}
Here $\varphi$ is the unique stationary solution with the same
mean value of the initial datum and $\beta$ is given by
\eqref{mejor.c.Poincare}.
\end{coro}

\begin{proof} It follows from
Theorem \ref{teo.Neumann.Lineal.conv.Exp} by considering that
$v=u-\varphi$ is a solution of \eqref{Neumann} with $h=0$.
\end{proof}


\section{Blow-up for $g(y,t) =
h(y)(T-t)^{-\alpha}$.}\label{SECT-b.u} \setcounter{equation}{0}

Now we analyze the asymptotic behavior of solutions of
\eqref{Neumann} when the flux at the boundary is given by
$$
g(y,t) = h(y) (T-t)^{-\alpha}
$$
with $h \geq 0$ and
$\int_\Omega\int_{\RR^N\setminus\Omega}J(x-y)h(y)\,dy>0$. We will
also assume that the initial data, and hence the solution, is non
negative. Throughout this section $u(x,t)$ will denote the
solution of \eqref{Neumann} with boundary and initial data as
described above. Also in this section we will assume, without loss
of generality, that $T < 1$. This makes the quantity $-\ln (T-t)$
positive which helps to avoid overloading the notation.

Throughout this section we will assume that $J\in L^\infty(\RR^N)$
and we will use the notation introduced in the Introduction.

First, we prove that $\alpha =1$ is the critical exponent to
obtain blowing up solutions.

\begin{lema}\label{explosion} The solution $u(x,t)$ blows up at
time $T$ if and only if $\alpha \ge 1$.
\end{lema}

\begin{proof} Set
$$
M(t) = \int_\Omega u(x,t) \, dx,
$$
then one has
$$
M' (t) = \frac{1}{(T-t)^\alpha}  \int_\Omega \int_{\RR^N \setminus
\Omega} J(x-y) h(y) \, dy \, dx \ge \frac{c}{(T-t)^\alpha}.
$$
Hence, if $\alpha \ge 1$ $M(t)$ is unbounded as $t\nearrow T$ and
the same is true for the solution $u(x,t)$.

On the other hand, if $\alpha <1$ we consider the solution of the
ordinary differential equation
$$
z' (t) = \frac{C}{(T-t)^\alpha} \mbox{ with } z(0) = z_0,
$$
that is a supersolution of our problem if $C$ and $z_0$ are large
enough. Since $z(t)$ remains bounded up to time $T$, a comparison
argument shows that so does $u(x,t)$.
\end{proof}

\begin{lema}\label{tasa por arriba} There exists a
constant $C$ such that for each integer $i$
such that $1 \leq i \leq \alpha$, the solution $u(x,t)$ verifies
$$
u(x,t) \leq \frac{C}{(T-t)^{\alpha -i }} \; \; \mbox{ in } \; \;
\Omega_{i-1} \; \mbox{ if } \; i \neq \alpha$$
and
$$
u(x,t) \leq -C\ln (T-t) \; \; \mbox{ in } \; \; \Omega_{i-1} \;
\mbox{ if } \; i = \alpha .$$
\end{lema}

\begin{proof} If $\alpha > 1$ we have that
$$
z(t) =  \frac{C_1}{(T-t)^{\alpha-1}}
$$
is a supersolution to our problem for $C_1$ large enough,
therefore
\begin{equation}\label{tasaB1.por.arriba}
u(x,t) \le \frac{C_1}{(T-t)^{\alpha-1}} \mbox{ in } \Omega_0.
\end{equation}

If $\alpha = 1$ the argument can be easily modified to get
$$
u(x,t) \leq -C_1\ln (T-t)  \; \; \mbox{ in } \; \; \Omega_{0}.$$

Now for $x \in \Omega_1 $ we have
\begin{equation} \label{Neumann.B2}
u_t (x,t) =   \displaystyle \int_{\Omega} J(x-y)(u(y,t) - u(x,t))
\,  dy
\end{equation}
which implies
\begin{equation} \label{Neumann.B2.2}
u_t (x,t) \leq  \displaystyle \int_{\Omega_1} J(x-y)( u
(y,t)-u(x,t)) \, dy \displaystyle + \int_{{\mathcal{B}}_1}
J(x-y) u(y,t) \, dy   \\
\end{equation}
Assume that $\alpha >2$. In this case, in view of
\eqref{tasaB1.por.arriba}, we can use the function
$$
z(t) = \frac{C_2 }{(T-t)^{\alpha-2}},
$$
with $C_2$ large enough, as a supersolution in $ \Omega_1$ to
obtain that
$$
u(x,t) \leq \frac{C_2}{(T-t)^{\alpha -2 }}, \qquad \mbox{ in }
\Omega_1.
$$
As before if $\alpha = 2$ we get
$$
u(x,t) \leq -C_2\ln (T-t),  \qquad \mbox{ in } \Omega_{1}.
$$

The previous argument can be repeated to obtain the conclusion of
the lemma with the constant $C = \max\limits_{1 \leq j \leq
[\alpha]}C_j$.
\end{proof}

We can describe now precisely the blow up set and profile of a
blowing up solution.

\begin{teo}
If $\alpha > 1$ is not an integer the blow-up set, $B(u)$, is
given by
$$
B(u) = \bigcup_{1\le i \leq [\alpha ]} {\mathcal B}_i,
$$
with the asymptotic behavior
$$
(T-t)^{\alpha -i } u(x,t) \to w_i (x) \; \mbox{ uniformly in }
{\mathcal B}_i \mbox{ as } \; t \to T$$ for each $i$ such that $1
\leq i \leq [\alpha ]$.

If $\alpha$ is an integer the blow-up set, $B(u)$, is given by
$$
B(u) = \bigcup_{1\le i \le \alpha} {\mathcal B}_i,
$$
with the asymptotic behavior,
$$
(T-t)^{\alpha -i } u(x,t) \to w_i (x) \; \mbox{ uniformly in }
{\mathcal B}_i \mbox{ as } \; t \to T$$ for each $i$ such that $1
\leq i < \alpha$ and
$$
\frac{u(x,t)}{-\ln (T-t)} \to \tilde w_\alpha (x) \; \mbox{
uniformly in } {\mathcal B}_\alpha \mbox{ as } \; t \to T.$$
\end{teo}

\begin{proof}

We have
\begin{equation} \label{Neumann.T-t}
u_t (x,t) =  \displaystyle \int_{\Omega} J(x-y)(u(y,t) - u(x,t))
\, dy + \displaystyle \int_{\RR^N \setminus \Omega} J(x-y)
\frac{h(y)}{(T-t)^\alpha} \, dy .
\end{equation}

We prove first the theorem in the case when $\alpha =1$.
Integrating \eqref{Neumann.T-t} in $t$ and using that, by Lemma
\eqref{tasa por arriba}, $u(x,t) \leq -C\ln (T-t)$ we get
$$\left|-\frac{u(x,t)}{\ln (T-t)} - \tilde w_1(x)\right|
\leq \frac{u(x,0)}{-\ln (T-t)} + C\frac{1}{-\ln
(T-t)}\displaystyle \int_0^t \ln (T-r) dr + \frac{\ln T}{-\ln
(T-t)}\tilde w_1(x).$$ This proves that
$$\lim\limits_{t \to T} \frac{u(x,t)}{-\ln (T-t)} = \tilde w_1(x)
\qquad \mbox{ uniformly in } \Omega_0.$$ Also if $x \in \Omega_1$
\eqref{Neumann.T-t} reads
$$u_t (x,t) =  \displaystyle \int_{\Omega} J(x-y)(u(y,t) - u(x,t))
\, dy .$$ Integrating in $t$ and using again Lemma \eqref{tasa por
arriba} we have
$$u(x,t) \leq u(x,0) - C\int_0^t \ln (T-r)dr.$$
Hence $u$ is bounded in $\Omega_1$ and the theorem is proved if
$\alpha =1$.

 Assume now that $\alpha > 1$ and consider the change of variables
$$
v_1(x,s) = (T-t)^{\alpha -1} u(x,t), \qquad s=-\ln (T-t).
$$
Since $u$ verifies \eqref{Neumann.T-t}, $v_1$ satisfies
\begin{equation} \label{Neumann.w}
\begin{array}{rl}
(v_1)_s (x,s) = & \displaystyle e^{-s} \int_{\Omega}
J(x-y)(v_1(y,s) - v_1(x,s)) \, dy
\\[12pt]
& + \displaystyle \int_{\RR^N \setminus \Omega} J(x-y) h(y) \, dy
 -(\alpha-1) v_1(x,s).
\end{array}
\end{equation}

Integrating in $s$ we obtain
\begin{equation} \label{Neumann.w.1}
\begin{array}{rl}
v_1 (x,s) -w_1(x) = & e^{-(\alpha -1)s}v(x,0)
\\[10pt]
 & + e^{-(\alpha -1)s} \displaystyle \int_0^s \displaystyle
e^{(\alpha -2)r} \int_{\Omega} J(x-y)(v_1(y,r) - v_1(x,r)) \,
dy\, dr \\[10pt]
& - e^{-(\alpha -1)s}w_1(x) .
\end{array}
\end{equation}

If $\alpha \neq 2$ since, by the previous lemma, $v_1$ is bounded
we get
\begin{equation} \label{Neumann.w.2}
|v_1 (x,s) -w_1(x)| \leq C( e^{-s}  + e^{-(\alpha -1)s})
\end{equation}
for some constant $C$. This implies that
\begin{equation}\label{asym1}
(T-t)^{\alpha -1 } u(x,t) \to w_1 (x)
\end{equation}
uniformly in $\Omega_{0}$ as $t \to T$.

We note that if $\alpha < 2$, since $w_1 (x)$ vanishes in
$\Omega_1$, \eqref{Neumann.w.2} implies
$$|v_1 (x,s)| \leq Ce^{-(\alpha -1)s} \; \mbox{ for } x \in \Omega_1$$
and hence
$$u(x,t) \leq C \; \mbox{ for } \; x \in \Omega_1.$$

Consequently, if $1< \alpha < 2$ the blow up set of $u$ is
$\Omega_0 \setminus \Omega_1 = {\mathcal B}_1$ and the asymptotic
behavior at the blow up time is given by \eqref{asym1}.

We have to handle now the case $\alpha = 2$ which is slightly
different. In this case instead of estimate \eqref{Neumann.w.2} there holds
\begin{equation} \label{Neumann.w.12}
|v_1 (x,s) -w_1(x)| \leq C( se^{-s}  + e^{-(\alpha -1)s}) .
\end{equation}
This still implies that
\begin{equation}\label{asym22}
(T-t)^{\alpha -1 } u(x,t) \to w_1 (x)
\end{equation}
uniformly in $\Omega_{0}$ as $t \to T$ but does not ensure that
$u$ is bounded in $\Omega_1$.

If $x \in \Omega_1$ one has
$$u_t (x,t) =  \displaystyle \int_{\Omega_1} J(x-y)(u(y,t) - u(x,t))
\, dy + \displaystyle \int_{{\mathcal B}_1} J(x-y) u(y,t) \, dy -
\displaystyle \int_{{\mathcal B}_1} J(x-y) u(x,t) \, dy .$$

Integrating in $t$ we obtain
$$u (x,t) - \int_0^t  \displaystyle \int_{{\mathcal B}_1} J(x-y) u(y,r) \,
dy\, dr = u(x,0) + I_1 + I_2,$$ where
$$I_1(x,t) = \int_0^t \displaystyle \int_{\Omega_1} J(x-y)(u(y,r) - u(x,r)) \,
dy\, dr
$$
and
$$I_2(x,t) = - \int_0^t \displaystyle \int_{{\mathcal B}_1} J(x-y) u(x,r) \,
dy\, dr.$$

Now using the fact that for $z \in \Omega_1$ one has $u(z,t) \leq
-C \ln (T-t)$, it can be checked that
$$\frac{I_1(x,t)}{\ln (T-t)} \to 0 \qquad \mbox{ uniformly in } \Omega_1
\mbox{ as } t \to T$$ and also
$$\frac{I_2(x,t)}{\ln (T-t)} \to 0 \qquad \mbox{ uniformly in } \Omega_1
\mbox{ as } t \to T.$$

Moreover since $(T-t) u(y,t) \to w_1(y)$ uniformly in ${\mathcal
B}_1$ as $t \to T$ one has
$$-\frac{1}{\ln (T-t)}\int_0^t  \displaystyle \int_{{\mathcal B}_1} J(x-y) u(y,r) \,
dydr \to \int_{{\mathcal B}_1} J(x-y) w_1(y) \, dy$$ uniformly in
$\Omega_1$ as $t \to T$.

Putting together this information we deduce that
$$-\frac{u(x,t)}{\ln (T-t)} \to \int_{{\mathcal B}_1} J(x-y) w_1(y) \, dy$$
uniformly in $\Omega_1$ as $t \to T$.

Finally since $u(x,t) \leq -C\ln (T-t)$ in $\Omega_1$ we can argue
as in the proof of the case $\alpha =1$ to show that $u$ remains
bounded in $\Omega_2$. So we have shown that if $\alpha = 2$, then
the blow up set, $B(u)$, of $u$ is given by
$$B(u) = {\mathcal B}_1 \cup {\mathcal B}_2$$
with the asymptotic behavior
$$\lim\limits_{t \to T}(T-t)u(x,t) = w_1(x) \qquad \mbox{ uniformly in } \; \Omega_0$$
and
$$\lim\limits_{t \to T}(T-t)u(x,t) = \tilde w_2(x) \qquad
\mbox{ uniformly in } \; \Omega_1.$$

If $\alpha > 2$ setting
$$
v_2(x,s) = (T-t)^{\alpha -2} u(x,t), \qquad s=-\ln (T-t),
$$
we obtain for $x \in \Omega_1$ the equation
\begin{equation} \label{Neumann.w.2.5}
(v_2)_s (x,s) = \displaystyle e^{-s} \int_{\Omega} J(x-y)(v_2(y,s)
- v_2(x,s)) \,  dy  -(\alpha-2) v_2(x,s).
\end{equation}

This can be written as
\begin{equation} \label{Neumann.w.3}
\begin{array}{rl}
(v_2)_s (x,s) = &  e^{-s} \displaystyle \int_{\Omega_1}
J(x-y)(v_2(y,s)-v_2(x,s)) \, dy \\[10pt]
& + \displaystyle \int_{B_1} J(x-y)v_1(y,s) \, dy -
v_1(x,s)\displaystyle \int_{B_1} J(x-y) \, dy
 \displaystyle
\\[10pt] &
-(\alpha-2) v_2(x,s).
\end{array}
\end{equation}
Again integrating in $s$ , after observing that by
\eqref{Neumann.w.2} \; $|v_1(x,s)| \leq Ce^{-s}$ since $x \in
\Omega_1$, we obtain that
\begin{equation} \label{Neumann.w.4}
\left|v_2 (x,s) -e^{-(\alpha - 2)s} \displaystyle \int_0^s
e^{-(\alpha - 2)r} \displaystyle \int_{{\mathcal B}_1}
J(x-y)v_1(y,r)\, dy\, dr \right| \leq C( e^{-s} + e^{-(\alpha
-2)s}) .
\end{equation}
for some constant $C$ provided that $\alpha \neq 3$.

Also by \eqref{Neumann.w.2} one has that
$$\left|\int_{{\mathcal B}_1} J(x-y)v_1(y,s)\, dy - \int_{{\mathcal B}_1}
J(x-y)w_1(y)\, dy \right| \leq Ce^{-s}$$ and consequently
$$\left|w_2(x) - e^{-(\alpha - 2)s} \displaystyle \int_0^s e^{-(\alpha
- 2)r} \displaystyle \int_{{\mathcal B}_1} J(x-y)v_1(y,r)\,dy \,
dr \right| \leq Ce^{-s}.$$ This, together with
\eqref{Neumann.w.4}, implies that for all $x \in \Omega_1$
$$ |v_2 (x,s) -w_2(x)| \leq C( e^{-s}  +  e^{-(\alpha
-2)s})
$$
for some constant $C$ and hence
$$
(T-t)^{\alpha -2 } u(x,t) \to w_2 (x)
$$
uniformly in $\Omega_{1}$.

The above procedure can be iterated to obtain for all integers $i$
such that $1 < i < \alpha$
\begin{equation}\label{vw}
|v_i (x,s) -w_i(x)|
\leq C( e^{-s} + e^{-(\alpha -i)s}) \; \; \mbox{ for all } \; \; x
\in \Omega_{i-1}
\end{equation}
for some constant $C$ and hence
\begin{equation}\label{asym2}
(T-t)^{\alpha -i } u(x,t) \to w_i (x)
\end{equation}
uniformly in $\Omega_{i-1}$. Moreover, it follows from \eqref{vw}
that for $x \in \Omega_{[\alpha ]}$, if $\alpha$ is not an
integer,
$$|v_{[\alpha]} (x,s)|
\leq Ce^{-(\alpha -[\alpha ])s}$$ and hence
$$u(x,t) \leq C \; \; \mbox{ for all } \; \; x \in
\Omega_{ [\alpha ]}.$$

In this fashion we have proved that, if $\alpha$ is not an
integer, the blow up set of $u$ is $\bigcup_{1 \leq i \leq [\alpha
]} {\mathcal B}_i$ and the behavior of $u$ near time $T$ in
${\mathcal B}_i$ is given by \eqref{asym2}. This proves the
theorem in the case that $\alpha$ is not an integer.

In the case that $\alpha$ is an integer one can argue as in the
proof of the case $\alpha = 2$ to obtain the result in that case.
\end{proof}


\section{Blow-up with a nonlinear boundary
condition.}\label{sect-NBC} \setcounter{equation}{0}

In this section we  deal with the problem
\begin{equation}\label{eqp}
\begin{array}{l}
\displaystyle u_t = \int_\Omega J(x-y)(u(y,t)-u(x,t))dy +
\int_{\Oc}
J(x-y)\overline{u}^p(y,t)dy.\\[12pt]
u(x,0) = u_0(x).
\end{array}
\end{equation}
Here we  assume that $J\in C(\RR^N)$, $u_0\in C(\overline\Omega)$, $u_0\ge0$ and
$\overline{u}$ is the extension of $u$ to a neighborhood of
$\overline{\Omega}$ defined as follows: take a small neighborhood
of $\partial \Omega$ in $\Oc$ in such a way that there exist
coordinates $(s, z) \in (0,s_0) \times
\partial\Omega$ that describe that neighborhood in the form $y=z +
s \eta (z)$ where $z\in \partial \Omega$ and $\eta(z)$ is the
exterior unit normal vector to $\partial \Omega$ at $z$. We set
$$
\overline{u} (y,t) = u(z,t).
$$
We also assume that $d < s_0$ therefore for any $x \in
\overline{\Omega}$ the ball centered at $x$ and of radius $d$ is
contained in the above mentioned neighborhood.

We address now the problem of local existence in time and
uniqueness of solutions.

As in the previous sections we set
$$A(x) = \int_\Omega J(x-y)dy$$
and observe that there exists $\alpha > 0$ such that $A(x) \geq
\alpha$ for all $x \in \overline{\Omega}$.

As earlier we  obtain a solution of \eqref{eqp} as a fixed
point of the operator $T$ defined by
$$
\begin{array}{rl}
Tu(x,t)= & \di e^{-A(x)t}u_0(x) + \int_0^t e
^{-A(x)(t-s)}\int_\Omega
J(x-y)u(y,s)\, dy\, ds \\[12pt]
& \di + \int_0^t e ^{-A(x)(t-s)}\int_{\Oc}
J(x-y)\overline{u}^p(y,s)\, dy\, ds.
\end{array}$$

We split the proof of existence into two cases. We deal first with
the case $p\ge 1$ since in this case we have uniqueness of
solutions. In this direction we have the following theorem.

\begin{teo}\label{locexistp>1}
$\mbox{ }$

\begin{itemize}
\item[a)] Let $p\ge 1$. There exists $t_0>0$ such that problem
\eqref{eqp} has a unique solution defined in $[0,t_0)$.

\smallskip

\item[b)] Let $p < 1$. There exists $t_0>0$ such that problem
\eqref{eqp} has at least one solution defined in $[0,t_0)$.
\end{itemize}
\end{teo}

\begin{proof}
Fix $M \geq ||u_0||_\infty$, $t_0 > 0$ and set $$X = \Big\{ u \in
C(\overline{\Omega} \times [0,t_0)) \; / \, u\ge0\;,\, |||u|||
\equiv \sup\limits_{(x,t) \in \overline{\Omega} \times [0,t_0) }
|u(x,t)| \leq 2M \Big\}.$$

If $t_0$ is chosen small enough, then $T$ maps $X$ into $X$.
Indeed, we have for $t \leq t_0$ and $u \in X$
$$
\begin{array}{rl}
|Tu(x,t)| &\leq  \displaystyle e^{-A(x)t}u_0(x) + \int_0^t e
^{-A(x)(t-s)}\int_\Omega J(x-y)|u(y,s)|\, dy\, ds \\[12pt]
& \displaystyle + \int_0^t e ^{-A(x)(t-s)}\int_{\Oc}
J(x-y)|\overline{u}^p(y,s)|\, dy\, ds \\[12pt]
& \leq M+ t_0(2M + (2M)^{p})\le 2M
\end{array}
$$
if $t_0$ is small.

Proof of a): We will prove that for $p \geq 1$ we can choose $t_0$
in such a way that $T$ is a strict contraction. In fact, for $t
\leq t_0$ and $u_1\; , \; u_2 \in X$
$$
\begin{array}{rl}
|Tu_1(x,t) - Tu_2(x,t)| & \leq \di \int_0^t e
^{-A(x)(t-s)}\int_\Omega
J(x-y)|u_1(y,s)-u_2(y,s)| \, dy\, ds \\[12pt]
& \di + \int_0^t e ^{-A(x)(t-s)}\int_{\Oc}
J(x-y)|\overline{u}_1^p(y,s)-\overline{u}_2^p(y,s)|\, dy\, ds
\\[12pt]
& \di \leq t_0(1 + p(2M)^{p-1})|||u_1-u_2|||
\end{array}
$$

\noindent and part a) of the theorem follows via Banach's fixed
point theorem.

\bigskip

Proof of b): We have that $T$ maps $X$ into $X$ if $t_0$ is small
enough. We claim that the operator $T:X \to X$ is compact. Indeed,
for $t_1, \; t_2 \leq t_0$,  $u \in X$ and $x_1, \; x_2 \in
\overline{\Omega}$ we have
$$
\begin{array}{l}
\displaystyle |Tu(x_1,t_1) - Tu(x_2,t_2)| \leq \Big| e
^{-A(x_1)t_1}u_0(x_1)-e ^{-A(x_2)t_2}u_0(x_2)\Big|\\
[12pt] +\Big|\displaystyle \int_0^{t_1}
e ^{-A(x_1)(t_1-s)}\int_\Omega J(x_1-y)u(y,s)\, dy\, ds \\[12pt]
\displaystyle \hskip0.5cm- \int_0^{t_2} e
^{-A(x_2)(t_2-s)}\int_\Omega J(x_2-y)u(y,s)\, dy\, ds \Big| \\[12pt]
\displaystyle \hskip1cm + \Big|\int_0^{t_1} e
^{-A(x_1)(t_1-s)}\int_{\Oc}
J(x_1-y)\overline{u}^p(y,s)\, dy\, ds \\[12pt]
\displaystyle \hskip1.5cm - \int_0^{t_2} e
^{-A(x_2)(t_2-s)}\int_{\Oc} J(x_2-y)\overline{u}^p(y,s)\, dy\, ds
\Big|.
\end{array}$$

As in the proof of Theorem \ref{comp.asimp} we have that for
$x_1,\; x_2 \in \overline{\Omega}$ one has
$$|e^{-A(x_1) t}-e^{-A(x_2) t}| \leq e^{-\alpha
t}t\ |A(x_1)-A(x_2)|$$ with $\alpha >0$. This inequality plus the
fact that $J$ is integrable imply, via a triangle inequality
argument, that the family $\{ Tu \; / \; u \in X \}$ is
equicontinuous and, since it is bounded, it is precompact in $(X,
||| \cdot |||)$. Consequently, since $T$ is clearly continuous in
$X$, it is a compact operator and the claim is proved. Part b) of
the theorem now follows from Schauder's fixed point theorem.
\end{proof}

\bigskip

\begin{rem} \label{remark.Lip} We observe that the same argument of
the proof of part {\em a)} of Theorem \ref{locexistp>1} provides
existence of a unique solution if the boundary nonlinearity takes
the form $f(\overline{u})$ with $f$ locally Lipschitz.
\end{rem}

Now we prove a comparison lemma for solutions of \eqref{eqp}.

\begin{lema}\label{comparisonp}
Let $u$ be a continuous subsolution and $v$ be a continuous
supersolution of problem \eqref{eqp} defined in $[0,t_0)$. Assume
$u(x,0) < v(x,0)$ for all $x \in \overline{\Omega}$. Then
$$u(x,t)
< v(x,t)$$ for all $(x,t) \in \overline{\Omega} \times [0,t_0).$
\end{lema}

\begin{proof}
Assume, for a contradiction, that the lemma is not true. Then, by
continuity, there exists $x_1 \in \overline{\Omega}$ and $0 < t_1
< t_0$ such that $u(x_1,t_1)=v(x_1,t_1)$ and $u(x,t) \leq v(x,t)$
for all $(x,t) \in \overline{\Omega} \times [0,t_1)$. We have now
$$
\begin{array}{rl}
0= & u(x_1,t_1)- v(x_1,t_1) = e
^{-A(x_1)t_1}(u(x_1,0)-v(x_1,0))
\\[8pt]
& \displaystyle +\int_0^{t_1} e ^{-A(x_1)(t_1-s)}\int_\Omega
J(x_1-y)\big(u(y,s)-v(y,s)\big)\, dy\, ds \\[12pt]
& \displaystyle + \int_0^{t_1} e ^{-A(x_1)(t_1-s)}\int_{\Oc}
J(x_1-y)\big(\overline{u}^p(y,s)-\overline{v}^p(y,s)\big)dyds < 0
\end{array}
$$ a
contradiction that proves the lemma.
\end{proof}

Now we use this comparison result to prove the lack of uniqueness
for $p<1$.

\begin{prop}\label{no.uni.p<1}
In the case $p<1$ with $u_0 \equiv 0$ there exists a nontrivial
solution of problem \eqref{eqp}. Hence this problem does not have
uniqueness.
\end{prop}

\begin{proof}
Let $b(t)$ be a positive solution of $b'= b^p$ with $b(0)=0$ and
$0\le a(x)\le\gamma$ be a continuous function with $a(x)\equiv\gamma$
on $\partial \Omega$.
Let $\gamma>0$ be so small as to have
$$
\gamma^p\int_{\RR^N\setminus\Omega}J(x-y)\,dy>2\gamma
$$
for every $x\in\Omega$.
Then,
$$
v(x,t) = a(x) b(t)
$$
is a subsolution to our problem for a certain interval of time,
$(0,t_0)$.

Let $\ve
>0$ be given and consider a locally Lipschitz function $f_\ve$
such that $f_\ve (s) = s^p$ for $s \geq \ve /2$. It follows from
Remark \ref{remark.Lip} that there  exists a unique solution,
$w_\ve$, of \eqref{eqp} with the boundary nonlinearity replaced by
$f_\ve (\overline{w})$ and initial data $w_\ve (x,0) \equiv \ve$.
By the comparison principle $w_\ve \geq \ve$ and hence it is a
solution of \eqref{eqp}.

By comparison, the sequence $w_\ve$ is monotone increasing in
$\ve$. In particular, for every $\ep$ $w_\ve$ is defined on the interval
$[0,t_1]$ where $w_1$ is. Therefore, by monotone
convergence, we obtain that the limit
$$
w= \lim_{\ve \to 0} w_\ve
$$
is a solution with $w(x,0)=0$.

 Using
again the comparison principle we obtain that $w_\ve (x,t) > v(x,t)$
for $0<t< \min \{t_0, t_1\}$. Hence, $w(x,t) >0$ for every
$0<t<\min \{t_0, t_1\}$ and all $x\in supp (a)$.
\end{proof}

We address now the blow up problem for solutions of \eqref{eqp}.
In this direction we have the following theorem.

\begin{teo}
$\mbox{ }$

\begin{itemize}
\item[a)] Let $p > 1$, then every non trivial
solution of \eqref{eqp} blows up in finite time.

\smallskip

\item[b)] Let $p \leq 1$, then every solution of \eqref{eqp} is
globally defined in time, by this we mean that it exists for all
$t \in [0,\infty )$.
\end{itemize}
\end{teo}

\begin{proof}

Proof of a): Let $u$ be a solution of \eqref{eqp} and assume, for
a contradiction, that it is globally defined in time.

Since $u \geq 0$ and $\int_\Omega J(x-y)\, dy \leq 1$ we have for
$x \in \overline{\Omega}$,
$$
u_t(x,t) \geq - u(x,t) + \int_{\Oc}J(x-y) \overline{u}^p(y,t)\,
dy.
$$
Here we have used that the equation is satisfied for $x \in
\partial \Omega$.

Integrating on $\partial \Omega$, denoting by $dS_x$ the surface
area element of $\partial \Omega$, we get
$$\frac{d}{dt}\int_{\partial \Omega}u(x,t)\, dS_x \geq -
\int_{\partial \Omega}u(x,t)\, dS_x +
\int_{\partial \Omega}\int_{\Oc}J(x-y) \overline{u}^p(y,t)\, dy\,
dS_x.$$

Since
$$\int_{\partial \Omega}\int_{\Oc}J(x-y) dydS_x
> 0$$
an application on Jensen's inequality implies that
$$\frac{d}{dt}\int_{\partial \Omega}u(x,t)\, dS_x \geq - \int_{\partial \Omega}u(x,t)\, dS_x +
C\left(\int_{\partial \Omega}\int_{\Oc}J(x-y) \overline{u}(y,t)\,
dy\, dS_x\right)^p$$ for some constant $C>0$.

Now,
$$
\begin{array}{l}
\displaystyle \int_{\partial \Omega}\int_{\Oc}J(x-y)
\overline{u}(y,t)\, dy \, dS_x \\[12pt]
\hskip0.5cm \displaystyle = \int_{\partial \Omega} \int_0^d
\int_{\partial \Omega} J\big(x- \sigma - s \eta (\sigma) \big)
u(\sigma, t) \, dS_\sigma \, ds \, dS_x \\[12pt]
\hskip1cm \displaystyle = \int_{\partial \Omega} \int_{\partial
\Omega} \Big[\int_0^d J\big(x- \sigma - s \eta (\sigma) \big) \,
ds \Big]
u(\sigma, t) \, dS_\sigma  \, dS_x \\[12pt]
\hskip1.5cm \displaystyle \ge \delta |\partial \Omega|
\int_{\partial \Omega} u(\sigma, t) \, dS_\sigma,
\end{array}
$$
since we have that $J(z) \ge c >0$ for $|z|<d/2$.

Thus, if we call
$$
m(t) = \int_{\partial \Omega} u(x,t) \, dS_x,
$$
we have
\begin{equation}\label{ec.bum}
m'(t) \ge - m(t) + \gamma m^p (t).
\end{equation}
This implies that $m(t) \to \infty$ in finite time if for some
$t_0$, $m(t_0)$ is large enough.

Since we are assuming that $u(x,t)$ is defined for every $t>0$, it
holds that $m(t)$ is defined (and finite) for all $t>0$. Let us
see that this leads to a contradiction. Let $v(x,t)$ be the
solution of \eqref{Neumann} with $g=0$ and $v(x,0)=u_0$. By
Theorem \ref{comp.asimp} we get that
$$
v(x,t) \to \displaystyle \frac{1}{|\Omega|} \int_\Omega u_0,
$$
uniformly in $\Omega$. Since $u$ is a supersolution for the
problem satisfied by $v$, there exists $t_1>0$ such that for $t
\ge t_1$,
$$
u(x,t) \ge \displaystyle \frac{1}{2|\Omega|} \int_\Omega u_0 = c_0
>0.
$$
Therefore,
$$
\begin{array}{rl}
M(t) = & \di \int_\Omega u(x,t) \, dx \\[12pt]
& \di \ge M(t_1) + \int_{t_1}^t \int_\Omega \int_{\Oc} J(x-y)
\overline{u}^p (y,s) \, dy \, dx \,
ds \\[14pt]
& \di \ge M(t_1)+ (t-t_1) c \ c_0^p.
\end{array}
$$

Arguing as before we get that
$$
u(x,t) \ge \frac{1}{2|\Omega|}M(t_2)
$$
for $t$ large enough. Now
$$
m(t_0)= \int_{\partial \Omega} u(x,t_0)\, dx \ge \frac{|\partial
\Omega|}{2|\Omega|}M(t_2) \ge C (M(t_1) + (t_2-t_1) c \ c_0^p).
$$
This implies that $m(t_0)$ is as large as we need if $t_0$ is
large enough, hence $m(t)$ is not defined for all times and we
conclude that $u$ blows up in finite time.

Proof of b): Let $p<1$ and let $u$ be a solution of \eqref{eqp}.
Set $v(t) = C(t+1)^{\frac{1}{1-p}}$. It is directly checked that
$$
v'(t) = \frac{C^{\frac{1}{1-p}}}{1-p}v^p(t).
$$
Picking $C$ such
that
$$
\frac{C^{\frac{1}{1-p}}}{1-p} \geq \max\limits_{x \in
\overline{\Omega}}\int_{\Oc} J(x-y)\, dy
$$
we have that $v$ is a
supersolution of \eqref{eqp}. Moreover taking $C$ larger, if
necessary, such that $u(x,0) < v(0)$ in $\overline{\Omega}$ we
obtain by Lemma \ref{comparisonp} that
$$u(x,t) \leq v(t)$$
as long as $u$ is defined. This implies the theorem in the case
$p<1$. The case $p=1$ is proved in the same fashion but using
$v(t) = Ce^t$ as a supersolution. \end{proof}

Our next result is an estimate of the blow up rate of blowing up
solutions of \eqref{eqp}.

\begin{teo}\label{teo.tasa.up}
Let $u$ be a solution of \eqref{eqp} that blows up at time $T$.
Then there exists a constant $C$ such that
\begin{equation}\label{tasa.up}
   (p-1)^{-1/(p-1)} (T-t)^{-1/(p-1)} \leq \|u (\cdot, t)
   \|_{L^\infty (\Omega)} \leq C (T-t)^{-1/(p-1)}.
\end{equation}
\end{teo}

\begin{proof}
Let
$$
v(t) = (p-1)^{-1/(p-1)} (T-t)^{-1/(p-1)}.
$$
One can easily check that $v$ is a supersolution of our problem.
If for some $t_0 \in (0,T)$ one has
$$
\|u (\cdot, t_0) \|_{L^\infty (\Omega)} < v(t_0),
$$
then there exists $\tilde{T} > T$ such that
$$
u(x,t_0) < (p-1)^{-1/(p-1)} (\tilde{T}-t_0)^{-1/(p-1)}.
$$
Let $\tilde v (t):=(p-1)^{-1/(p-1)} (\tilde{T}-t)^{-1/(p-1)}$
that is also a supersolution
to our problem in the interval $[t_0 , \tilde T)$.

Using a comparison argument, we obtain
$$
u(x,t) < (p-1)^{-1/(p-1)} (\tilde{T}-t)^{-1/(p-1)},
$$
for all $t\in (t_0,T)$. This contradicts the fact that $\tilde{T} > T$ and
hence
$$
  (p-1)^{-1/(p-1)}
(T-t)^{-1/(p-1)} \leq \| u(\cdot, t) \|_{L^\infty (\Omega)}.
$$

The proof of the reverse inequality is more involved. By equation
\eqref{ec.bum}, if
$$
m(t) = \int_{\partial \Omega} u(x,t) \, dS_x \to \infty \qquad
\mbox{ as } t \nearrow T,
$$
we have
\begin{equation}\label{bum.m.tasa}
m(t) \le C(T-t)^{-1/(p-1)}.
\end{equation}

Now, we claim that
\begin{equation}\label{Cosa.acotada}
    (T-t)^{1/(p-1)} \int_0^t \int_{\Oc} J(x-y) \overline{u}^p
    (y,s) \, dy \, ds \le C,
\end{equation}
for all $x \in \partial \Omega$. In fact, if this does not hold,
there exists a sequence $(x_n,t_n)$ with $x_n \in \partial
\Omega$, $t_n \nearrow T$, such that
$$
(T-t_n)^{1/(p-1)} \int_0^{t_n} \int_{\Oc} J(x_n-y) \overline{u}^p
    (y,s) \, dy \, ds \to \infty.
$$
By compactness we may assume that $x_n \to x_0 \in \partial
\Omega$. Hence
$$
(T-t_n)^{1/(p-1)} \int_0^{t_n} \int_{(\Oc ) \cap B(x_0, 2d)}
\overline{u}^p
    (y,s) \, dy \, ds \to \infty.
$$
Therefore there exists a point $x_1 \in \partial \Omega$ such that
for a subsequence that we still call $t_n$,
$$
(T-t_n)^{1/(p-1)} \int_0^{t_n} \int_{(\Oc ) \cap B(x_1, d/4)}
\overline{u}^p
    (y,s) \, dy \, ds \to \infty.
$$
Since every function involved is nonnegative and $J(z) \ge c
>0$ for $|z|<d/2 $ we get
$$
(T-t_n)^{1/(p-1)} \int_0^{t_n} \int_{\Oc} J(\hat{x}-y)
\overline{u}^p
    (y,s) \, dy \, ds \to \infty,
$$
for every $\hat{x}\in \partial \Omega \cap \{|\hat{x}-x_1|< d
/4\}$.

Using \eqref{repre.u}, we get
$$
(T- t_n )^{1/(p-1)} u(\hat{x}, t_n) \ge c (T-t_n)^{1/(p-1)}
\int_0^{t_n} \int_{\Oc} J(\hat{x}-y) \overline{u}^p (y,s) \, dy \,
ds \to \infty.
$$
Therefore,
$$
(T-t_n)^{1/(p-1)}m(t_n) = (T-t_n)^{1/(p-1)} \int_{\partial \Omega}
u(x,t_n) \, dS_x \to \infty,
$$
which contradicts \eqref{bum.m.tasa}. The claim is proved.

Using again that $J(z) \ge c
>0$ for $z<d/2 $ we get that \eqref{Cosa.acotada} holds for every $x \in
\overline{\Omega}$. In fact, first we see that for every $x\in\partial\Omega$
$$
(T-t)^{1/(p-1)}\int_0^t\int_{\partial\Omega\cap B_{d/4}(x)}
u^p(\sigma,s)\,dS_\sigma\,ds\le C.
$$
Then, since $\partial \Omega$ is compact we deduce that
$$
(T-t)^{1/(p-1)}\int_0^t\int_{\partial\Omega}u^p(\sigma,s)\,dS_\sigma\,ds\le C.
$$
This immediately implies, by using that $J\in L^\infty$, that
\eqref{Cosa.acotada} holds for every $x \in
\overline{\Omega}$.

Now, let for $t_0 <T$,
$$
M = \max_{\overline{\Omega} \times [0,t_0]} (T-t)^{1/(p-1)}u(x,t)
= (T-t_1)^{1/(p-1)}u(x_1,t_1).
$$
This implies by using again \eqref{repre.u} that
$$
M \le C + \int_0^{t_1} e^{-A(x_1) (t_1-s)} \int_\Omega J(x_1 -y )
M \, dy \, ds \le C + (1- e^{-A(x_1) t_1}) M.
$$
So that, since $A(x) \ge \alpha >0$,
$$
M \le C,
$$
with $C$ independent of $t_0$. The result follows.
\end{proof}

\begin{coro}\label{conj.explosion}
Let $u$ be a solution of \eqref{eqp} that blows up at time $T$.
Then, the blow-up set, $B(u)$, verifies
\begin{equation}\label{conj.up}
B(u) \subset \{ x \in \Omega \; / \; \mbox{dist }(x, \partial
\Omega ) \leq Kd \}
\end{equation}
where $K = [p/(p-1)]$.
\end{coro}

\begin{proof}
The proof follows from the results in Section \ref{SECT-b.u}.
\end{proof}

{\bf Acknowledgements.} Part of this work was done during visits
of JDR and NW to Universidad Cat\'olica de Chile. These authors are
grateful for the warm hospitality.

\end{document}